\documentstyle[amscd,amssymb,verbatim,11pt]{amsart}

\theoremstyle{plain}
\newtheorem{Prop}{Proposition}[section]
\newtheorem{Thm}[Prop]{Theorem}
\newtheorem{Cor}[Prop]{Corollary}

\newtheorem{Lem}[Prop]{Lemma}

\theoremstyle{definition}
\newtheorem{Def}[Prop]{Definition}

\theoremstyle{remark}

\newtheorem{Problem}[Prop]{\bf Problem}

\def\dim{\mathop{\roman{dim}}}
\def\int{\mathop{\roman{int}}}

\def\1{^{-1}}

\def\dim{\text{dim}}

\def\mesh{\text{mesh}}

\def\dist{\text{dist}}
\def\NN{{\mathcal N}}
\def\UU{{\mathcal U}}
\def\VV{{\mathcal V}}
\def\WW{{\mathcal W}}
\def\RR{{\mathbb R}}
\def\dokaz{{\bf Proof. }}
\numberwithin{equation}{section}
\begin{document}
\title[
Nagata-Assouad dimension via Lipschitz extensions
]%
   {Nagata-Assouad dimension via Lipschitz extensions
}

\author{N.~Brodskiy}
\address{University of Tennessee, Knoxville, TN 37996, USA}
\email{brodskiy@@math.utk.edu}

\author{J.~Dydak}
\address{University of Tennessee, Knoxville, TN 37996, USA}
\email{dydak@@math.utk.edu}

\author{J.~Higes}
\address{Departamento de Geometr\'{\i}a y Topolog\'{\i}a,
Facultad de CC.Matem\'aticas.
Universidad Complutense de Madrid.
Madrid, 28040 Spain}
\email{josemhiges@@yahoo.es}

\author{A.~Mitra}
\address{University of Tennessee, Knoxville, TN 37996, USA}
\email{ajmitra@@math.utk.edu}

\date{ January 10, 2006}


\thanks{ The second-named author was partially supported
by Grant No.2004047  from the United States-Israel Binational Science
Foundation (BSF),  Jerusalem, Israel.
}

\begin{abstract} 
In the first part of the paper we show how to relate several dimension theories (asymptotic dimension with Higson property, asymptotic dimension of Gromov, and capacity dimension of Buyalo \cite{Buyalo1}) to Nagata-Assouad dimension. This is done by applying two functors on the Lipschitz category of metric spaces: microscopic and macroscopic. In the second part we identify (among spaces of finite Nagata-Assouad dimension) spaces of Nagata-Assouad dimension at most $n$ as those for which the $n$-sphere $S^n$ is a Lipschitz extensor.
Large scale and small scale analogs of that result are given.
\end{abstract}

\maketitle

\medskip
\medskip
\tableofcontents

\section{Introduction}\label{section Introduction}

The large-scale geometry of metric spaces has been the subject of intense
research in the last 15 years. For a comprehensive account of the area see
Gromov's paper \cite{Gro asym invar}. There are two large-scale concepts related to the topic
of our paper: {\it asymptotic dimension} and {\it asymptotic dimension of linear type}
(Dranishnikov and Zarichnyi  \cite{Dran-Zar} refer to the latter as {\it asymptotic dimension with Higson
property}). Asymptotic dimension is an invariant of the coarse category of Roe \cite{Roe lectures}.
Asymptotic dimension of linear type is preserved by bi-Lipschitz functions, so its natural place is in the Lipschitz category.
\par Notice that our notion of a Lipschitz function is a bit different from that used in
\cite{HeinonenBook}. Namely, we allow Lipschitz constant to be smaller than $1$. Thus,
a Lipschitz function $f:(X,d_X)\to (Y,d_Y)$ satisfies $d_Y(f(x),f(y))\leq \lambda\cdot d_X(x,y)$
for some $\lambda\ge 0$ and all $x,y\in X$. The infimum of all possible $\lambda$ is denoted by $Lip(f)$. $f$ is called {\it bi-Lipschitz} if there are constants $\mu,\lambda > 0$
such that $\mu\cdot d_X(x,y)\leq d_Y(f(x),f(y))\leq \lambda\cdot d_X(x,y)$
for all $x,y\in X$.

\par
In \cite{Lang-Sch Nagata dim} a variation 
of asymptotic dimension is considered. That invariant of bi-Lipschitz functions was introduced and named {\it Nagata dimension} by
Assouad \cite{Assouad Sur la distance de Nagata} as it is closely related to a theorem of Nagata characterizing
the topological dimension of metrizable spaces (cf. [\cite{Nagata Note on dimension theory}, Thm. 5] or [\cite{Nagata Modern Dimension Theory},
p. 138]). In our paper we refer to it as {\it Nagata-Assouad dimension}.
In contrast to the asymptotic dimension, the Nagata dimension
of a metric space is in general not preserved under quasi-isometries, but
it is still a bi-Lipschitz invariant and, as it turns out, even a quasisymmetry
invariant (see \cite{Lang-Sch Nagata dim}). The class of metric spaces with finite Nagata dimension
includes all doubling spaces, metric trees, euclidean buildings, and homogeneous
or pinched negatively curved Hadamard manifolds as shown in \cite{Lang-Sch Nagata dim}.
One of main results of \cite{Lang-Sch Nagata dim} relates to theorems of Assouad \cite{Assouad Plongements lipschitziens} and
Dranishnikov \cite{Dran hypersphericity}. Namely, Theorem 1.3 of \cite{Lang-Sch Nagata dim} states: Every metric space with Nagata-Assouad dimension at most $n$ admits a quasisymmetric embedding into the product of $n + 1$
metric trees. 
\par
The result from \cite{Lang-Sch Nagata dim} of major importance to us is Theorem 1.4:
\begin{Thm} \label{LangSThm}
Suppose that $X$, $Y$ are metric spaces, $\dim_{NA} X \leq n < \infty$,
and $Y$ is complete. If $Y$ is Lipschitz $m$-connected for $m = 0, 1,\ldots , n-1$,
then the pair $(X,Y)$ has the Lipschitz extension property.
\end{Thm}

Recall that $(X,Y)$ has the {\it Lipschitz extension property} if there is a constant $C > 0$
such that for any Lipschitz map $f:A\to Y$, $A$ any subset of $X$, there is a Lipschitz
extension $g:X\to Y$ of $f$ such that $Lip(g)\leq C\cdot Lip(f)$.
We call $Y$ a {\it Lipschitz extensor of $X$} in such a case.
$Y$ is {\it Lipschitz $m$-connected} if there is a constant $C_m > 0$ such that
any Lipschitz function $f:S^m\to Y$ extends over the $(m+1)$-ball $B^{m+1}$
to $g:B^{m+1}\to Y$ so that $Lip(g)\leq C_m\cdot Lip(f)$.

One of the main themes of our paper is characterizing Nagata-Assouad dimension
via Lipschitz extensions. In Section~\ref{section Lipschitz extensions and Nagata-Assouad dimension}
we characterize spaces of Nagata-Assouad dimension at most $n$ (among all spaces
of finite Nagata-Assouad dimension) as those for which $n$-sphere $S^n$ is
a Lipschitz extensor. In the case of dimension zero the assumption
of $\dim_{NA}(X)$ being finite can be dropped.

Recently, S.Buyalo \cite{Buyalo1} introduced the {\it capacity dimension} of a metric space
and proved many analogs of results obtained by
 U. Lang and T. Schlichenmaier \cite{Lang-Sch Nagata dim} for Nagata-Assouad
 dimension. It is clear that capacity dimension is the small scale version of Nagata-Assouad dimension.
 In Section~\ref{section Microscopic and macroscopic Nagata-Assouad dimensions}
 we formalize that observation by introducing microscopic and macroscopic functors on
 the Lipschitz category. That way many results from Section 3 of \cite{Buyalo1}
 can be deduced formally from \cite{Lang-Sch Nagata dim}.
 Also, the main result of \cite{Buyalo1} (the asymptotic dimension of a visual hyperbolic space $X$
 is bounded by 1 plus the capacity dimension of the visual boundary of $X$) is really about
 Nagata-Assouad dimension of the visual boundary of $X$ - see \ref{MicroNagataForBoundedIsNagata}.


\section{Microscopic and macroscopic Nagata-Assouad dimensions}\label{section Microscopic and macroscopic Nagata-Assouad dimensions}
\begin{Def} \label{DefOfNagataDim}
A metric space $X$ is said to be of {\it Nagata-Assouad dimension}
at most $n$ (notation: $\dim_{NA}(X)\leq n$) if there is $C > 0$ such that for all $r > 0$
there is a cover $\UU=\bigcup\limits_{i=1}^{n+1}\UU_i$ of $X$
so that each $\UU_i$ is $r$-disjoint and the diameter of elements of $\UU$ is
bounded by $C\cdot r$.
\end{Def}

Nagata-Assouad dimension can be characterized in many ways (see \cite{Lang-Sch Nagata dim}):
\begin{Prop} \label{PossibleCharOfNagata}
For a metric space $(X,d)$ the following conditions are equivalent:
\begin{enumerate}
\item $\dim_{NA}(X)\leq n$.
\item There is a constant $C_1 > 0$ such that for any $r > 0$
there is a cover $\UU_r$ of $X$ of multiplicity at most $n+1$,
of mesh at most $C_1\cdot r$, and of Lebesque number at least $r$.
\item There is a constant $C_2 > 0$ such that for any $r > 0$
there is a cover $\VV_r$ of mesh at most $C_2\cdot r$
such that each $r$-ball $B(x,r)$ intersects at most $n+1$ elements of
$\VV_r$.
\end{enumerate}
\end{Prop}

If one replaces all $r > 0$ in \ref{DefOfNagataDim} by $r$ sufficiently small,
then one gets the concept of {\it capacity dimension} of Buyalo \cite{Buyalo1}
for which he proved analog of \ref{PossibleCharOfNagata}.
If one replaces all $r > 0$ in \ref{DefOfNagataDim} by $r$ sufficiently large,
then one gets the concept of {\it asymptotic dimension of linear type}
or {\it asymptotic dimension with Higson property}.
In this section we will show how to introduce those dimensions formally
from Nagata-Assouad dimension. That way small scale version and large scale version
of \ref{PossibleCharOfNagata} are in fact consequences of \ref{PossibleCharOfNagata}.

Given a metric space $(X,d)$ and $\epsilon > 0$ we consider the metric $\max(d,\epsilon)$
on $X$. Needles to say, the formula should not be read literally, only in the case of $x\ne y$.
Similarly, we consider the metric $\min(d,\epsilon)$.

A metric space $(X,d)$ is called {\it discrete} if $(X,d)$ is {\it $\delta$-discrete} for some $\delta > 0$. That means $d(x,y) > \delta$ for all $x\ne y$.

\begin{Lem} \label{BiLipEquivalenceForMacroFunctor}
Any discrete metric space $(X,d)$ is bi-Lipschitz equivalent to $(X,\max(d,\epsilon))$
for all $\epsilon > 0$.
\end{Lem}
\dokaz Suppose $(X,d)$ is $\delta$-discrete. Notice the identity map
\par\noindent $id:(X,\max(d,\epsilon))\to (X,d)$ is $1$-Lipschitz and its inverse is $(1+\frac{\epsilon}{\delta})$-Lipschitz as
$d(x,y)\leq \max(d(x,y),\epsilon))\leq (1+\frac{\epsilon}{\delta})\cdot d(x,y)$ for all $x\ne y\in X$.
\hfill $\blacksquare$

\begin{Cor} \label{MacroDoesNotDependOnEpsilon}
For any metric space $(X,d)$  and $\epsilon,\delta > 0$ the space $(X,\max(d,\epsilon))$ is
bi-Lipschitz equivalent to $(X,\max(d,\delta))$.
\end{Cor}
\dokaz Assume $\epsilon > \delta$ and notice $\max(d_\delta,\epsilon)=d_\epsilon$,
where $d_a=\max(d,a)$. Use \ref{BiLipEquivalenceForMacroFunctor}.
\hfill $\blacksquare$

Since Nagata-Assouad dimension is an invariant of the Lipschitz category
(see \cite{Lang-Sch Nagata dim} for a stronger result for quasisymmetric embeddings),
one gets the following.
\begin{Cor} \label{NagataDimDoesNotDependOnEpsilon}
For any metric space $(X,d)$ the Nagata-Assouad dimension of $(X,\max(d,\epsilon))$
does not depend on $\epsilon > 0$.
\end{Cor}

Given a metric space $(X,d)$ one can disregard its microscopic features
by considering the space $(X,\max(d,1))$.
The {\it macroscopic Nagata-Assouad dimension} of $(X,d)$ is defined
as $\dim_{NA}(X,\max(d,1))$.

\begin{Cor} \label{NagataEqualsMacroForDiscrete}
If $(X,d)$ is a discrete metric space, then its macroscopic Nagata-Assouad dimension
equals the Nagata-Assouad dimension $\dim_{NA}(X,d)$ of $(X,d)$.
\end{Cor}

\begin{Lem} \label{CharOfMacroNagata}
The macroscopic Nagata-Assouad dimension of a metric space $X$ is at most $n\ge 0$ if and only if there is $C > 0$ such that for sufficiently large $r > 0$
there is a cover $\UU=\bigcup\limits_{i=1}^{n+1}\UU_i$ of $X$
so that each $\UU_i$ is $r$-disjoint and the diameter of elements of $\UU$ is
bounded by $C\cdot r$.
\end{Lem}
\dokaz Suppose $X$ has a constant $C > 0$ such that for all $r > M$, where $M > 0$,
there is a cover $\UU_r=\bigcup\limits_{i=1}^{n+1}\UU_i$ of $X$
so that each $\UU_i$ is $r$-disjoint and the diameter of elements of $\UU_r$ is
bounded by $C\cdot r$.
Let $d_M=\max(d,M)$. Notice the cover $\UU_r$, $r\leq M$, consisting of one-point sets
is $r$-disjoint in $(X,d_M)$. Since covers $\UU_r$, $r > M$, have the same desired properties
in $(X,d_M)$ as they do in $(X,d)$, $\dim_{NA}(X,d_M)\leq n$.
\par If $\dim_{NA}(X,d_1)\leq n$, then for each $r> 1$ there is a a cover $\UU_r=\bigcup\limits_{i=1}^{n+1}\UU_i$ of $X$
so that each $\UU_i$ is $r$-disjoint (in $(X,d_1)$) and the diameter of elements of $\UU_r$ is
bounded by $C\cdot r$ in $(X,d_1)$.
Notice that $\UU_i$ is also $r$-disjoint in $(X,d)$ and the diameter of elements of $\UU_r$ is
bounded by $(C+1)\cdot r$ in $(X,d)$.
\hfill $\blacksquare$

In view of definition of the Higson property in \cite{Dran asym top} one has the following consequence of \ref{CharOfMacroNagata}.
\begin{Cor} \label{MacroNagataEqualsHigson}
If $(X,d)$ is a metric space, then the macroscopic Nagata-Assouad dimension of $(X,d)$ is at most $n$ if and only if
$asdim(X)\leq n$ with the Higson property.
\end{Cor}

In the reminder of this section we will dualize the above results from large scale/macroscopic category
to small scale/microscopic category.

\begin{Lem} \label{BoundedSpacesEquivalentToMicro}
Any bounded metric space $(X,d)$ is bi-Lipschitz equivalent to $(X,\min(d,\epsilon))$
for all $\epsilon > 0$.
\end{Lem}
\dokaz Suppose $(X,d)$ is $\delta$-bounded. Notice the identity map
$id:(X,\min(d,\epsilon))\to (X,d)$ is $(1+\frac{\epsilon}{\delta})$-Lipschitz  and its inverse is $1$-Lipschitz as
$\min(d(x,y),\epsilon)\leq d(x,y)\leq (1+\frac{\delta}{\epsilon})\cdot \min(d(x,y),\epsilon))$ for all $x\ne y\in X$.
\hfill $\blacksquare$

\begin{Cor} \label{MicroDoesNotDependOnEpsilon}
For any metric space $(X,d)$  and $\epsilon,\delta > 0$ the space $(X,\min(d,\epsilon))$ is
bi-Lipschitz equivalent to $(X,\min(d,\delta))$.
\end{Cor}
\dokaz Assume $\epsilon < \delta$ and notice $\min(d_\delta,\epsilon)=d_\epsilon$,
where $d_a=\min(d,a)$. Use \ref{BoundedSpacesEquivalentToMicro}.
\hfill $\blacksquare$

Since Nagata-Assouad dimension is an invariant of the Lipschitz category,
one gets the following.

\begin{Cor} \label{NagataDimIndependentOnMicro}
For any metric space $(X,d)$ the Nagata-Assouad dimension of $(X,\min(d,\epsilon))$
does not depend on $\epsilon > 0$.
\end{Cor}

Given a metric space $(X,d)$ one can disregard its macroscopic features
by considering the space $(X,\min(d,1))$.
The {\it microscopic Nagata-Assouad dimension} of $(X,d)$ is defined
as $\dim_{NA}(X,\min(d,1))$.

\begin{Cor} \label{MicroNagataForBoundedIsNagata}
If $(X,d)$ is a bounded metric space, then its microscopic Nagata-Assouad dimension
equals the Nagata-Assouad dimension $\dim_{NA}(X,d)$ of $(X,d)$.
\end{Cor}

\begin{Lem} \label{CharOfMicroNagata}
The microscopic Nagata-Assouad dimension of a metric space $X$ is at most $n$ if and only if there is $C > 0$ such that for sufficiently small $r > 0$
there is a cover $\UU=\bigcup\limits_{i=1}^{n+1}\UU_i$ of $X$
so that each $\UU_i$ is $r$-disjoint and the diameter of elements of $\UU$ is
bounded by $C\cdot r$.
\end{Lem}
\dokaz Suppose $X$ has a constant $C > 0$ such that for all $r < M$, where $M > 0$,
there is a cover $\UU_r=\bigcup\limits_{i=1}^{n+1}\UU_i$ of $X$
so that each $\UU_i$ is $r$-disjoint and the diameter of elements of $\UU_r$ is
bounded by $C\cdot r$.
Let $d_M=\min(d,M)$. Notice the cover $\UU_r$, $r\ge M$, consisting of the whole $X$
has diameter at most $(C+1)\cdot r$ in $(X,d_M)$. Since covers $\UU_r$, $r < M$, have the same desired properties
in $(X,d_M)$ as they do in $(X,d)$, $\dim_{NA}(X,d_M)\leq n$.
\par If $\dim_{NA}(X,d_1)\leq n$, then for each $r < 1$ there is a a cover $\UU_r=\bigcup\limits_{i=1}^{n+1}\UU_i$ of $X$
so that each $\UU_i$ is $r$-disjoint (in $(X,d_1)$) and the diameter of elements of $\UU_r$ is
bounded by $C\cdot r$ in $(X,d_1)$.
Notice that $\UU_i$ is also $r$-disjoint in $(X,d)$ and the diameter of elements of $\UU_r$ is
bounded by $C\cdot r$ in $(X,d)$.
\hfill $\blacksquare$

Since the capacity dimension of Buyalo \cite{Buyalo1} can be characterized
by the condition appearing in \ref{CharOfMicroNagata}, one derives
the following.
\begin{Cor} \label{MicroEqualsCapacity}
If $(X,d)$ is a metric space, then the microscopic Nagata-Assoud dimension of $(X,d)$ is at most $n$ if and only if
the capacity dimension of $X$ is at most $n$.
\end{Cor}

In \cite{Lang-Sch Nagata dim} it is shown that if $X=A\cup B$, then $$\dim_{NA}(X)=\max(\dim_{NA}(A),\dim_{NA}(B)).$$

\begin{Cor} \label{UnionThmForCaporHigs}
Let $D(Y)$ stand for either the capacity dimension of $Y$ or
the asymptotic dimension of linear type of $Y$.
If $X=A\cup B$, then $$D(X)=\max(D(A),D(B)).$$
\end{Cor}
\dokaz Switch to either $\max(d,1)$ or $\min(d,1)$ metrics. 
\hfill $\blacksquare$

\section{Spheres as Lipschitz extensors}\label{section Spheres as Lipschitz extensors}

A metric space $E$ is a {\it Lipschitz extensor}
of $X$ if there is a constant $C > 0$ such that any $\lambda$-Lipschitz function
$f:A\to E$, $A$ a subset of $X$, extends to a $C\cdot\lambda$-Lipschitz function $\tilde f:X\to E$.

The purpose of this section is to find necessary and sufficient conditions
for a sphere $S^m$ to be a Lipschitz extensor of $X$. This is done by comparing existence of Lipschitz extensions in a finite range of Lipschitz
constants to existence of Lebesque refinements in a finite
range of Lebesque constants (see \ref{CharOfSmAsLipExtensor}, \ref{SmBeingLipExtensorInARangeImpliesCovers}, and \ref{nCoversTonPlus1Covers}).

Given a cover $\UU=\{U_s\}_{s\in S}$ of a metric space $(X,d)$
there is a natural family of functions $\{f_s\}_{s\in S}$ associated to $\UU$:
$f_s(x):=\dist(x,X\setminus U_s)$. To simplify matters by {\it the local Lebesque number} $L_\UU(x)$
of $\UU$ at $x$ we mean $$\sup\{f_s(x)\mid s\in S\}$$ and by the (global) {\it Lebesque number}
$L(\UU)$ of $\UU$ we mean $$\inf\{L_\UU(x)\mid x\in X\}.$$
We are interested in covers with positive Lebesque number. For those
{\it the local multiplicity} $m_\UU(x)$ can be defined as  
$1+|T(x)|$, where $T(x)=\{s\in S\mid f_s(x) > 0\}$
and the global multiplicity $m(\UU)$ can be defined as $$\sup\{m_\UU(x)\mid x\in X\}.$$
If the multiplicity $m(\UU)$ is finite, then $\UU$ has a natural partition of unity
$\{\phi_s\}_{s\in S}$ associated to it:
$$\phi_s(x)=\frac{f_s(x)}{\sum\limits_{t\in S} f_t(x)}.$$
That partition can be considered as a {\it barycentric map} $\phi:X\to\NN(\UU)$
from $X$ to the {\it nerve} of $\UU$.
Since each $f_s$ is $1$-Lipschitz, $\sum\limits_{t\in S} f_t(x)$ is $2m(\UU)$-Lipschitz
and each $\phi_s$ is $\frac{2m(\UU)}{L(\UU)}$-Lipschitz (use the fact that
$\frac{u}{u+v}$ is $\frac{\max(Lip(u),Lip(v))}{\inf(u+v)}$-Lipschitz). Therefore 
$\phi:X\to\NN(\UU)$ is $\frac{4m(\UU)^2}{L(\UU)}$-Lipschitz. See \cite{BD} and \cite{BS} for more details and better estimates of Lipschitz constants.

Most estimates in this paper work well for both the $l_1$ and $l_2$ metrics on $R^n$ and simplicial complexes.

\par In analogy to $\lambda$-Lipschitz functions we introduce the concept of
{\it $r$-Lebesque cover} $\UU$. That is simply a shortcut to $r\leq L(\UU)$.

 \begin{Prop} \label{CharOfSmAsLipExtensor}
Suppose $X$ is metric,
 $m\ge 0$, $C > 0$, and $\lambda_2 > \lambda_1 > 0$.
If any $\lambda$-Lipschitz function
$f:A\to S^m$, $A$ a subset of $X$ and $\lambda_1 < \lambda < \lambda_2$, extends to a $C\cdot\lambda$-Lipschitz function $\tilde f:X\to S^m$, then $t=\frac{1}{4C(m+2)^2(m+1)}$ has the property that any finite 
$r$-Lebesque cover $\UU=\{U_0,\ldots,U_{m+1}\}$ of $X$
admits a refinement $\VV$ so that $\VV$ is $t\cdot r$-Lebesque and the multiplicity
of $\VV$ is at most $m+1$ provided $\frac{4(m+2)^2}{\lambda_2} < r < \frac{4(m+2)^2}{\lambda_1}$.
\end{Prop}
\dokaz Assume $\frac{4(m+2)^2}{\lambda_2} < r < \frac{4(m+2)^2}{\lambda_1}$
and $\UU$ is $r$-Lebesque. 
Therefore $\lambda_1 < \frac{4(m+2)^2}{r} < \lambda_2$.
Consider a barycentric map $\phi:X\to \NN(\UU)=\Delta^{m+1}$.
Notice $Lip(\phi)\leq \frac{4(m+2)^2}{r}$.
There is $g:X\to \partial\Delta^{m+1}$ such that $Lip(g)\leq \frac{4C(m+2)^2}{r}$ and
$g(x)=\phi(x)$ for all $x\in X$ so that $\phi(x)\in \partial\Delta^{m+1}$. 
Consider $V_i=\{x\in X\mid g_i(x) > 0\}$. 
Notice $\VV=\{V_i\}_{i=0}^{i=m+1}$ is of multiplicity at most $m+1$.
Also $x\in V_i$ implies $x\in U_i$, so $\VV$ refines $\UU$.
Given $x\in X$ there is $i$ such that $g_i(x)\ge \frac{1}{m+1}$.
If $d(x,y) < \frac{r}{4C(m+2)^2(m+1)}$, then
$|g_i(x)-g_i(y)| < \frac{1}{m+1}$ and $g_i(y) > 0$. Thus, the ball at $x$ of radius
$\frac{r}{4C(m+2)^2(m+1)}$ is contained in one element of $\VV$.
That proves $\VV$ is $t\cdot r$-Lebesque, where $t=\frac{1}{4C(m+2)^2(m+1)}$.
\hfill $\blacksquare$

In this paper we consider the space $\RR^n$ with either the
$l_1$-metric $d_1(x,y)=\sum\limits_{i=1}^n|x_i-y_i|$ or the $l_2$-metric
$d_2(x,y)=\sqrt{\sum\limits_{i=1}^n(x_i-y_i)^2}$. The inequalities
$$d_1(x,y)\le n\cdot\max\{|x_i-y_i|\}\le n\cdot d_2(x,y)$$
$$d_2(x,y)\le \sqrt{n}\cdot\max\{|x_i-y_i|\}\le \sqrt{n}\cdot d_1(x,y)$$
show that the identity map $(\RR^n,d_2)\to(\RR^n,d_1)$ is
$n$-Lipschitz and its inverse is $\sqrt{n}$-Lipschitz.

\begin{Lem} \label{simplex is absolute Lip extensor}
Let $\triangle$ be a closed convex subset of the space
$(\RR^n,d_j)$, $j=1\text{ or }2$. For any metric space $X$ any
$\lambda$-Lipschitz map $f\colon A\to \triangle$ of a subspace
$A\subset X$ can be extended to a $n^2\cdot\lambda$-Lipschitz map
$\tilde{f}\colon X\to \triangle$.
\end{Lem}
\dokaz Fix an orthogonal coordinate system in $(\RR^n,d_2)$. Given
a $\lambda$-Lipschitz map $f\colon A\to (\RR^n,d_2)$, every
coordinate map $f_i\colon A\to\RR$ is $\lambda$-Lipschitz and can
be extended to a $\lambda$-Lipschitz map $\hat{f_i}\colon
X\to\RR$~\cite{McShane}. These coordinate extensions define the
map $\hat{f}\colon X\to\RR^n$ which is
$\lambda\sqrt{n}$-Lipschitz. Clearly, the nearest point retraction
$r_\triangle\colon (\RR^n,d_2)\to\triangle$ is 1-Lipschitz.
Therefore the composition $\tilde{f}=r_\triangle\circ\hat{f}$ is
$\sqrt{n}\lambda$-Lipschitz.

Since the nearest point retraction $r_\triangle\colon
(\RR^n,d_1)\to\triangle$ may be multivalued, we proceed as
follows. The composition $id\circ f\colon A\rightarrow
(\triangle,d_1)\to (\triangle,d_2)$ is $\sqrt{n}\lambda$-Lipschitz
and admits a ${n}\lambda$-Lipschitz extension $\tilde{f}_2\colon
X\to (\triangle,d_2)$ by the first part of the proof. Then the
composition $\tilde{f}=id\circ \tilde{f}_2\colon X\to
(\triangle,d_2)\to(\triangle,d_1)$ is $n^2\lambda$-Lipschitz.
\hfill $\blacksquare$

The idea behind the proof of the next proposition is best understood if
one thinks of maps from $X$ to an $(m+1)$-simplex $\Delta^{m+1}$
as a partition of unity. Since we want to create a map to its boundary
$S^m=\partial \Delta^{m+1}$, a geometrical tool is the radial projection $r$
which we splice in the form of $(1-\beta)\cdot r+ \beta\cdot \phi$
with a partition of unity $\phi$ coming from a covering of $X$
of multiplicity at most $m+1$.

 \begin{Prop} \label{SmBeingLipExtensorInARangeImpliesCovers}
Suppose $X$ is metric, $m\ge 0$, $t > 0$, and $r_2 > r_1 > 0$. There is $s > 0$ such that
if any finite $r$-Lebesque cover $\UU=\{U_0,\ldots,U_{m+1}\}$ of $X$, where $r_1 < r < r_2$,
admits an $t\cdot r$-Lebesque refinement $\VV$  satisfying $m(\VV)\leq m+1$,
then any $\lambda$-Lipschitz function $f:A\to S^m$, $A\subset X$,
admits a $C\cdot\lambda$-Lipschitz extension $\tilde f:X\to S^m$
provided $\frac{1}{12sr_2 (m+2)} < \lambda < \frac{1}{12sr_1(m+2)}$
and $C=50(m+2)^2s+\frac{150s^2 (m+2)^5}{t}$.
\end{Prop}
\dokaz
$s > 0$ is chosen so that given a $\lambda$-Lipschitz $f:A\to \Delta^{m+1}$
one can
extend it to an $s\cdot\lambda$-Lipschitz $g:X\to \Delta^{m+1}$
(see \ref{simplex is absolute Lip extensor}). 

\par Suppose $\frac{1}{12sr_2 (m+2)} < \lambda < \frac{1}{12sr_1(m+2)}$ and
$f:A\to \partial\Delta^{m+1}$ is $\lambda$-Lipschitz.
Extend it to an $s\cdot\lambda$-Lipschitz $g:X\to \Delta^{m+1}$.
Let $\alpha:X\to [0,1]$ be defined as $\alpha(x)=(m+2)\cdot \min\{g_i(x)\mid 0\leq i\leq m+1\}$.
Notice $Lip(\alpha)\leq (m+2)s\cdot \lambda$. 
Let $\beta:[0,1]\to [0,1]$ be defined by $\beta(z)=3z-1$ on $[1/3,2/3]$,
$\beta(z)=0$ for $z\leq 1/3$ and $\beta(z)=1$ for $z\ge 2/3$.
Notice $Lip(\beta)\leq 3$.

Put $U_i=\{x\in X\mid g_i(x) > \frac{\alpha(x)}{m+2} \text{ or } \alpha(x)> 2/3\}$
and notice $L(\UU)\ge r= \frac{1}{12s\lambda (m+2)}$ as follows:
\par Case 1: $x\in X$ and $\alpha(x) > 3/4$. Now, for any $y\in X$
with $d(x,y) < \frac{1}{12s\lambda (m+2)}$ one has
$\alpha(x)-\alpha(y) \leq 1/12$, so $\alpha(y) \leq 2/3$ is not possible.
Thus, in that case the ball $B(x,\frac{1}{12s\lambda (m+2)})$ is contained in all $U_i$.
\par
Case 2: $\alpha(x)\leq 3/4$. There is $i$ so that $g_i(x)\ge \frac{1}{m+2}$.
Since $\psi_i=g_i-\frac{\alpha}{m+2}$ is $2s\lambda$-Lipschitz,
for any $y\in X$ satisfying $d(x,y) < \frac{1}{8s\lambda (m+2)}$
one has $\psi_i(x)-\psi_i(y)< \frac{1}{4(m+2)}$
and $\psi_i(y) > 0$ as $\psi_i(x) \ge \frac{1}{4(m+2)}$.
\par
Thus $\UU$ is $r$-Lebesque and $r_1 < r < r_2$.
Shrink each $U_i$ to $V_i$ so that $m(\VV)\leq m+1$
and $L(\VV)\ge \frac{t}{12s\lambda (m+2)}$.
The barycentric map $\phi:X\to \partial\Delta^{m+1}$
corresponding to $\VV$ has $Lip(\phi)\leq \frac{4(m+2)^2}{L(\VV)}\leq \frac{48s\lambda (m+2)^3}{t}$.

Define $h(x)=\sum\limits_{i=0}^{m+1} (g_i(x)-\frac{\alpha(x)}{m+2})\cdot \frac{1-\beta(\alpha(x))}{1-\alpha(x)}\cdot e_i+\sum\limits_{i=0}^{m+1} \beta(\alpha(x))\cdot \phi_i(x)\cdot e_i$.
To show $Lip(h)\leq C\cdot\lambda$
we will use the following observations. 
\begin{enumerate}
\item If $u,v:X\to [0,M]$, then $Lip(u\cdot v)\leq M\cdot(Lip(u)+Lip(v))$.
\item In addition, if $v:X\to [k,M]$ and $k > 0$, then $$Lip(\frac{u}{v})\leq M\cdot \frac{Lip(u)+Lip(v)}{k^2}.$$
\item $v(x)=1-\alpha(x)\ge 1/3$ if $\frac{1-\beta(\alpha(x))}{1-\alpha(x)} > 0$.
\end{enumerate}

Therefore $Lip(\sum\limits_{i=0}^{m+1} \beta(\alpha(x))\cdot \phi_i(x)\cdot e_i)\leq 3\cdot Lip(\alpha)\cdot (m+2)\cdot Lip(\phi)\leq 3(m+2)\cdot (m+2)s\cdot \lambda\cdot\frac{48s\lambda (m+2)^3}{t}\leq
\frac{150s^2\lambda (m+2)^5}{t}$.
Also, $Lip(\frac{1-\beta(\alpha(x))}{1-\alpha(x)})\leq 9\cdot 4\cdot Lip(\alpha)\leq 36(m+2)s\lambda$,
so $Lip(\sum\limits_{i=0}^{m+1} (g_i(x)-\frac{\alpha(x)}{m+2})\cdot \frac{1-\beta(\alpha(x))}{1-\alpha(x)})\leq
(m+2)\cdot (2s\lambda+36(m+2)s\lambda)\leq 50(m+2)^2s\lambda$
and $C=50(m+2)^2s+\frac{150s^2 (m+2)^5}{t}$ works.
\par
It remains to show $h(X)\subset \partial\Delta^{m+1}$ and $h|A=f$.
$h|A=f$ follows from the fact $\alpha(x)=0$ if $x\in A$.
It is clear $h(x)\in \partial\Delta^{m+1}$ if either $\beta(\alpha(x))=0$ or $\beta(\alpha(x))=1$,
so assume $0 < \beta(\alpha(x)) < 1$. 
In that case $\phi_i(x) > 0$ implies $g_i(x)-\frac{\alpha(x)}{m+2} > 0$,
so the only possibility for $h(x)$ to miss $\partial\Delta^{m+1}$
is when $g_i(x)-\frac{\alpha(x)}{m+2} > 0$ for all $i$ which is not possible.
\hfill $\blacksquare$

Propositions \ref{CharOfSmAsLipExtensor} and \ref{SmBeingLipExtensorImpliesCovers} imply the following

 \begin{Cor} \label{SmBeingLipExtensorImpliesCovers}
If $X$ is metric and
 $m\ge 0$, then the following conditions are equivalent:
 \begin{itemize}
\item[a.] $S^m$ is a Lipschitz extensor of $X$.
\item[b.] There is $t > 0$ such that any finite cover $\UU=\{U_0,\ldots,U_{m+1}\}$ of $X$
admits a refinement $\VV$ so that $L(\VV)\ge tL(\UU)$ and the multiplicity
of $\VV$ is at most $m+1$.
\end{itemize}
\end{Cor}

\begin{Prop} \label{nCoversTonPlus1Covers}
Suppose $X$ is a metric space, $n\ge 0$, $1> t > 0$, and $r_2 > r_1 > 0$. If 
every $r$-Lebesque cover $\UU=\{U_0,\ldots,U_{n+1}\}$ of $X$, where $r_1 < r < r_2$,
admits a $4t\cdot r$-Lebesque refinement $\VV$ satisfying $m(\VV)\leq n+1$,
then any $s$-Lebesque cover $\WW=\{W_0,\ldots, W_{n+2}\}$ of $X$, where $4r_1 < s < 4r_2$,
admits a $t\cdot s$-Lebesque refinement $\VV$ of multiplicity at most $n+2$.
\end{Prop}
\dokaz Suppose $4r_1 < s < 4r_2$. First, let us show that any $s/2$-Lebesque cover
$\UU=\{U_i\}_{i=0}^{i=n+1}$ of $A\subset X$ consisting of $n+2$ elements has a refinement $\VV$ such that
$L(\VV)\ge t\cdot s$ and $m(\VV)\leq n+1$.
Define $U'_i=U_i\cup (X\setminus A)$ for $i\leq n+1$. Notice $L(\UU')\ge s/4$.
Indeed, if $x\in X$, then $B(x,s/4)\cap A$ is either empty
or is contained in $B(y,s/2)$ for some $y\in A$.
Since $B(y,s/2)\cap A\subset U_i$ for some $i\leq n+1$,
$B(y,s/2)\subset U_i'$ and $B(x,s/4)\subset U_i'$.
There is a cover $\WW$ of $X$ such that $\WW$ refines $\UU'$, $L(\WW)\ge 4t\cdot s/4$,
and $m(\WW)\leq n+1$. Putting $\VV=\WW|_A$ finishes the task.
\par
Suppose $\WW=\{W_0,\ldots,W_{n+2}\}$ is an $s$-Lebesque cover of $X$.
Let $A$ be the union of balls $B(x,s/2)$ such that $B(x,s)$ is not contained in $W_{n+2}$.
Define $U_i=W_i\cap A$ for $i\leq n+1$ and observe $L(\UU)\ge s/2$
for $\UU=\{U_i\}_{i=0}^{i=n+1}$ as a cover of $A$. Indeed, if $x\in A$,
then there is $y\in X$ such that $B(y,s)$ is not contained in $W_{n+2}$
and $x\in B(y,s/2)$. Therefore, $B(y,s)\subset W_i$ for some $i\leq n+1$
which means $B(x,s/2)\cap A\subset B(y,s)\cap A\subset W_i\cap A=U_i$.
\par
Shrink each $U_i$ to $V_i$ so that the intersection of all $V_i$ is empty and
$L(\VV)\ge t\cdot s$. 
Define $W'_i=V_i$ for $i\leq n+1$ and $W'_{n+2}=W_{n+2}$.
The cover $\WW'$ is of multiplicity at most $n+2$.
We want to show $L(\WW')\ge t\cdot s$.
If $B(x,s)\subset W_{n+2}$, we are done. Otherwise $B(x,s/2)\subset A$ and
there is $i\leq n+1$ such that $B(x,t\cdot s)\subset V_i$ in which case $B(x,t\cdot s)\subset W'_i$.
\hfill $\blacksquare$

\begin{Cor} \label{HigherSpheresAreLipAE}
Suppose $X$ is a metric space and $n\ge 0$. If $S^{n}$ is a Lipschitz extensor of $X$, then so is $S^{n+1}$.
\end{Cor}
\dokaz By \ref{SmBeingLipExtensorImpliesCovers} there is $t > 0$ such that any cover
$\UU$ of $X$ consisting of $n+2$ elements has a refinement $\VV$ satisfying
$L(\VV)\ge t\cdot L(\UU)$ and $m(\VV)\leq n+1$. 
Use $\lambda_2$ very large and $\lambda_1$ very small.
We may assume $t < 1$, so applying \ref{nCoversTonPlus1Covers}
and \ref{SmBeingLipExtensorImpliesCovers} completes the proof.
\hfill $\blacksquare$

\begin{Def}\label{LargeSmallScaleLipExtensors}
A metric space $E$ is a {\it large scale Lipschitz extensor}
(respectively, {\it a small scale Lipschitz extensor})
of $X$ if there are constants $C, M > 0$ such that any $\lambda$-Lipschitz function
$f:A\to E$, $A$ a subset of $X$, extends to a $C\cdot\lambda$-Lipschitz function $\tilde f:X\to E$ for all $\lambda < M$
(respectively, for all $\lambda > M$).
\end{Def}

Using \ref{CharOfSmAsLipExtensor}, \ref{SmBeingLipExtensorInARangeImpliesCovers}, and \ref{nCoversTonPlus1Covers} one proves easily the following
large/small scale analogs of \ref{SmBeingLipExtensorImpliesCovers}.
 \begin{Cor} \label{SmBeingScaleLipExtensorEquivToCovers}
If $X$ is metric and
 $m\ge 0$, then the following conditions are equivalent:
 \begin{itemize}
\item[a.] $S^m$ is a large scale Lipschitz extensor (respectively, a small scale Lipschitz extensor) of $X$.
\item[b.] There are constants $t, M > 0$ such that any finite 
$r$-Lebesque cover $\UU=\{U_0,\ldots,U_{m+1}\}$ of $X$,
where $r > M$ (respectively, $r < M$),
admits a $t\cdot r$-Lebesque refinement $\VV$ of multiplicity
at most $m+1$.
\end{itemize}
\end{Cor}

As in \ref{HigherSpheresAreLipAE} one proves its large/small scale analogs:
\begin{Cor} \label{HigherSpheresAreScaleLipAE}
Suppose $X$ is a metric space and $n\ge 0$. If $S^{n}$ is a large scale Lipschitz extensor (respectively, a small scale Lipschitz extensor) of $X$, then so is $S^{n+1}$.
\end{Cor}

\section{Lipschitz extensions and Nagata-Assouad dimension}\label{section Lipschitz extensions and Nagata-Assouad dimension}

\begin{Thm} \label{CharNADimAmongFiniteDimensional}
Suppose $X$ is a metric space of finite Nagata-Assouad dimension.
If 
 $n\ge 0$, then the following conditions are equivalent:
 \begin{itemize}
\item[a.] $S^n$ is a Lipschitz extensor of $X$.
\item[b.] $\dim_{NA}(X)\leq n$.
\end{itemize}
\end{Thm}
\dokaz The direction (b)$\implies$(a) follows from  \ref{SmBeingLipExtensorImpliesCovers} as follows.
Given a cover $\UU$ of $X$ of Lebesque number $L(\UU) > 0$
pick $\VV$ of $\mesh(\VV) < L(\UU)$ so that $L(\VV)> L(\UU)/C$ and $m(\VV)\leq m+1$.
Notice that $\VV$ refines $\UU$ and $L(\VV)> L(\UU)/C$.
\par
(a)$\implies$(b).
Without loss of generality (in view of \ref{HigherSpheresAreLipAE}), we may assume $\dim_{NA}(X)\leq n+1$. $k > 0$ is chosen so that given a $\lambda$-Lipschitz $f:A\to \partial\Delta^{n+1}$
one can
extend it to an $k\cdot\lambda$-Lipschitz $g:X\to \partial\Delta^{n+1}$. 
Let $c > 1$ be a constant such that for any $r > 0$ there is a cover $\UU$ of $X$
of mesh at most $c\cdot r$, Lebesque number at least $r$, that can be expressed as a union $\bigcup\limits_{i=1}^{n+2}\UU_i$
so that each $\UU_i$ is $r$-disjoint.
For such a cover pick a barycentric map $f:X\to \mathcal{N}(\UU)$ to the nerve
of $\UU$ such that $Lip(f)\leq \frac{4(n+2)^2}{r}$.
Given an $(n+1)$-simplex $\Delta$ in $\mathcal{N}(\UU)$ we look at
$f|f^{-1}(\partial \Delta)$ and extend it over $f^{-1}(\Delta)$ to obtain $g_\Delta:f^{-1}(\Delta)\to \partial\Delta$ of Lipschitz number at most $k\cdot \frac{4(n+2)^2}{r}$. Paste all $g_\Delta$ together
to $g:X\to \mathcal{N}(\UU)$. Our goal it to estimate the mesh and Lebesque number
of $g^{-1}(st(v))$, $v$ a vertex of $\mathcal{N}(\UU)$. 

\par The mesh of $\{g^{-1}(st(v))\}$ is at most twice that of $\{f^{-1}(st(v))\}$. Indeed, if
$g(x)\in st(w)$ and $g(x)\ne f(x)$, then $f(x)$ must belong to the interior of a simplex $\Delta$
containing $w$. Thus $g(x)$ belongs to the star of $st(w)$ in $\{st(v)\}$ and $x$ belongs
to the star of $U_w$ in $\UU$ ($U_w$ is the element of $\UU$ corresponding to $w$).
That star is of size at most $4cr$. 
\par Suppose $C$ is a subset of $X$ of diameter less than $\frac{r}{4k\cdot (n+3)^3}$. 
Pick all elements
$U_0,\ldots, U_m$ of $\UU$ intersecting $C$ and let $v_i$ be corresponding vertices
of $\mathcal{N}(\UU)$. Notice $m\leq n+1$ and $g(C)$ is contained in the simplex $[v_0,\ldots,v_m]$
of $\mathcal{N}(\UU)$. Pick $x_0\in C$ and, without loss of generality, assume the barycentric coordinate $\phi_{v_0}(x_0)$ of $g(x_0)$ corresponding to $v_0$ is at least $\frac{1}{n+2}$.
Suppose $g(x)$ does not belong to $st(v_0)$ for some $x\in C$. Thus $\phi_{v_0}(x)=0$ and
$\frac{1}{n+2}\leq |\phi_{v_0}(x)-\phi_{v_0}(x_0)|\leq (k\cdot \frac{(n+3)^2}{r})\cdot d(x,x_0)
\leq (k\cdot \frac{4(n+3)^2}{r})\cdot \frac{r}{4k\cdot (n+3)^3}=\frac{1}{n+3}$, a contradiction.

\par By setting $d=\frac{r}{4k\cdot (n+3)^3}$ the argument above shows the existence
of a cover $\VV$ of $X$ of multiplicity at most $n+1$, of Lebesque number at least $d$, and of mesh at most $d\cdot 16c\cdot k\cdot (n+3)^3$. That means $\dim_{NA}(X)\leq n$.
\hfill $\blacksquare$

Adjusting the proof of \ref{CharNADimAmongFiniteDimensional}
one can deduce the following.
\begin{Thm} \label{CharScaleNADimAmongFiniteDimensional}
Suppose $X$ is a metric space of finite capacity dimension
(respectively, of finite asymptotic dimension with Higson property).
If 
 $n\ge 0$, then the following conditions are equivalent:
 \begin{itemize}
\item[a.] $S^n$ is a small scale (respectively, large scale) Lipschitz extensor of $X$.
\item[b.] The capacity dimension
(respectively, the asymptotic dimension with Higson property) of $X$
is at most $n$.
\end{itemize}
\end{Thm}

\begin{Problem} \label{CharNAviaSnQuestion}
Suppose $X$ is a metric space such that $S^n$ is a Lipschitz extensor of $X$.
Is $\dim_{NA}(X)$ at most $n$?
\end{Problem}


In \cite{Brod-Dydak-Higes-Mitra} the authors proved the following.

\begin{Thm} \label{NADimViaExtendingToSZero}
For a metric space $X$ the following conditions are equivalent:
\begin{itemize}
\item[a.]  $\dim_{NA}(X)\leq 0$.
\item[b.] Every metric space $Y$ is a Lipschitz extensor of $X$.
\item[c.] The $0$-sphere $S^0$ is a Lipschitz extensor of $X$.\end{itemize}
\end{Thm}
Thus, for $n=0$,
the answer to \ref{CharNAviaSnQuestion} is positive.

In the remainder of this section we extend \ref{NADimViaExtendingToSZero}
to large and small scales.
\begin{Cor} \label{CapacityDimZeroChar}
For a metric space $X$ the following conditions are equivalent:
\begin{itemize}
\item[a.]  The capacity dimension of $X$ is at most $0$.
\item[b.] Each bounded metric space $Y$ is a small scale Lipschitz extensor of $X$.
\item[c.] The $0$-sphere $S^0$ is a small scale Lipschitz extensor of $X$.
\end{itemize}
\end{Cor}
\dokaz a)$\implies$b). Pick $\epsilon > 0$.
Let $d=\min(d_X,\epsilon)$ and let $Y$ be $S$-bounded.
Define the relation between $M$ and $\epsilon$ by $M=\frac{S}{\epsilon}$
and notice that for $\lambda\ge M=\frac{S}{\epsilon}$, a function
$f:(A,d_X|A)\to (Y,d_Y)$ is $\lambda$-Lipschitz if and only if
the function $f:(A,d|A)\to Y$ is $\lambda$-Lipschitz.
\par
Let $C > 1$ be a constant such that any $\lambda$-Lipschitz $f:(A,d|A)\to (Y,d_Y)$
extends to $C\cdot\lambda$-Lipschitz $\tilde f:(X,d)\to (Y,d_Y)$.
Given $\lambda$-Lipschitz function
$f:(A,d_X|A)\to (Y,d_Y)$, $\lambda\ge M=\frac{S}{\epsilon}$,
the function $f:(A,d|A)\to Y$ is also $\lambda$-Lipschitz.
Extend it to a $C\cdot\lambda$-Lipschitz $\tilde f:(X,d)\to Y$
and notice that $\tilde f:(X,d_X)\to Y$ is also $C\cdot\lambda$-Lipschitz.
\par b)$\implies$c) is obvious.
\par c)$\implies$a). Let $M > 1$ be a number such that any $f:A\to S^0$ satisfying
$Lip(f)\ge M$ has an extension $g:X\to S^0$ so that $Lip(g)\leq C\cdot Lip(f)$, where $C > 1$.
Suppose $r< \frac{1}{(C+1)M}$. Consider the equivalence classes determined
by $x\sim y$ if and only if $x$ can be connected to $y$ by a chain of points separated by at most
$r$. If any of them has diameter bigger that $Cr$, then there are points $x$ and $y$ in that particular
class such that $Cr\leq d(x,y)\leq (C+1)r$. Pick injection $f:\{x,y\}\to S^0$. Its Lipschitz constant is at least $\frac{1}{Cr+r} >  M$ and at most $\frac{1}{Cr}$.
By extending it to $g:B\to S^0$ of Lipschitz constant at most $\frac{1}{r}$ we arrive at a contradiction
that points of the chain joining $x$ and $y$ are mapped to the same point by $g$. By  \ref{CharOfMicroNagata} the capacity dimension of $X$ is at most $0$.
\hfill $\blacksquare$

\begin{Cor} \label{AsympDimHigsonZeroChar}
For a metric space $X$ the following conditions are equivalent:
\begin{itemize}
\item[a.]  The asymptotic dimension of $X$ is at most $0$ with the Higson property.
\item[b.] Each discrete metric space $Y$ is a large scale Lipschitz extensor of $X$.
\item[c.] The $0$-sphere $S^0$ is a large scale Lipschitz extensor of $X$.
\end{itemize}
\end{Cor}
\dokaz a)$\implies$b). Pick $\epsilon > 0$.
Let $d=\max(d_X,\epsilon)$ and let $Y$ be $\delta$-discrete.
Define the relation between $M$ and $\epsilon$ by $M=\frac{\delta}{\epsilon}$
and notice that for $\lambda\leq M=\frac{\delta}{\epsilon}$, a function
$f:(A,d_X|A)\to (Y,d_Y)$ is $\lambda$-Lipschitz if and only if
the function $f:(A,d|A)\to Y$ is $\lambda$-Lipschitz.
\par
Let $C > 1$ be a constant such that any $\lambda$-Lipschitz $f:(A,d|A)\to (Y,d_Y)$
extends to $C\cdot\lambda$-Lipschitz $\tilde f:(X,d)\to (Y,d_Y)$.
Given $\lambda$-Lipschitz function
$f:(A,d_X|A)\to (Y,d_Y)$, $\lambda\leq M=\frac{\delta}{C\epsilon}$,
the function $f:(A,d|A)\to Y$ is also $\lambda$-Lipschitz.
Extend it to a $C\cdot\lambda$-Lipschitz $\tilde f:(X,d)\to Y$
and notice that $\tilde f:(X,d_X)\to Y$ is also $C\cdot\lambda$-Lipschitz.
\par b)$\implies$c) is obvious.
\par c)$\implies$a). Let $M > 0$ be a number such that any $f:A\to S^0$ satisfying
$Lip(f)\leq M$ has an extension $g:X\to S^0$ so that $Lip(g)\leq C\cdot Lip(f)$, where $C > 1$.
Suppose $r> \frac{1}{CM}$. Consider the equivalence classes determined
by $x\sim y$ if and only if $x$ can be connected to $y$ by a chain of points separated by less than
$r$. If any of them has diameter bigger that $Cr$, then there are points $x$ and $y$ in that particular
class such that $Cr > d(x,y)$. Pick injection $f:\{x,y\}\to S^0$. Its Lipschitz constant is less than $\frac{1}{Cr} <  M$.
By extending it to $g:B\to S^0$ of Lipschitz constant less than $\frac{1}{r}$ we arrive at a contradiction
that points of the chain joining $x$ and $y$ are mapped to the same point by $g$. By  \ref{MacroNagataEqualsHigson} the asymptotic dimension of $X$ with Higson property is at most $0$.
\hfill $\blacksquare$

A way to probe solving \ref{CharNAviaSnQuestion} would be to investigate,
for a given $n > 0$, the class of metric spaces $X$ such that $S^n$ is a Lipschitz extensor of $X$.
One faces immediately the question of extending \ref{UnionThmForCaporHigs}:
\begin{Problem} \label{UnionProblemForLipExtensors}
Suppose $X=A\cup B$ is a metric space such that $S^n$ is a Lipschitz extensor of $A$ and $B$.
Is $S^n$ is a Lipschitz extensor of $X$?
\end{Problem}

\section{Coarsely equivalent metrics and Nagata-Assouad dimension}\label{section Coarsely equivalent metrics and Nagata-Assouad dimension}

In this section we characterize asymptotic dimension of Gromov in terms of
Nagata-Assouad dimension.
\begin{Thm} \label{CharOfAsdimViaNagata}
For an unbounded metric space $(X,d)$ the following conditions are equivalent:
\begin{itemize}
\item[a.]  $asdim(X)\leq n$.

\item[b.] There is a hyperbolic metric $(X,d_h)$ coarsely equivalent to $(X,d)$
such that $\dim_{NA}(X,d_h)\leq n$ and the Gromov boundary $\partial_\infty X$
of $X$ consists of one point.

\item[c.] There is a metric $(X,d_1)$ coarsely equivalent to $(X,d)$
such that Nagata-Assouad dimension $\dim_{NA}(X,d_1)$ of $(X,d_1)$ is at most $n$.

\end{itemize}
\end{Thm}
\dokaz a)$\implies$b). Pick a sequence of covers $\UU_i$ of $X$, $i\ge 1$,
of multiplicity at most $n+1$ such that $\mesh(\UU_i)\to\infty$, $L(\UU_i)\to \infty$,
and $2\mesh(\UU_i) <  L(\UU_{i+1})$ for all $i$.
If $x\ne y$, define $d_h(x,y)$ as the smallest $i$ so that there is $U\in\UU_i$ containing both $x$ and $y$.
Clearly, $d_h$ is coarsely equivalent to $d$: if $d_h(x,y)\leq i$, then $d(x,y)\leq \mesh(\UU_i)$.
Also, $d(x,y)\leq L(\UU_i)$ implies $d_h(x,y)\leq i$.
 \par
 Notice that for any triangle in $(X,d_h)$ with sides $a\ge b\ge c$ one has $a\leq b+1$.
 The reason for this is that $x,y\in U\in\UU_i$ and $y,z\in V\in\UU_i$ implies existence
 of $W\in\UU_{i+1}$ containing all three points $x,y,z$ as $U\cup V$ is of diameter less than
 the Lebesque number of $\UU_{i+1}$.
 Thus, in any triangle of $(X,d_h)$ the difference of any two sides that are not minimal is either $-1$, $0$,
 or $1$.
\par Fix $x_0\in X$ and consider the Gromov product 
$$(x|y)=\frac{d_h(x,x_0)+d_h(y,x_0)-d_h(x,y)}{2}.$$
To show $(X,d_h)$ is Gromov hyperbolic it suffices to prove
$$(x|z)\ge \min((x|y),(y|z))-4$$
 for all $x,y,z\in X$. Equivalently, the smallest product in a triangle is at least the medium one minus $4$.
\par If all distances $d_h(x,x_0)$, $d_h(y,x_0)$, and $d_h(z,x_0)$ are within $1.5$ from a number $t$,
then, as $d_h(z,x)$ cannot be larger than both $d_h(y,x)+1$ and $d_h(y,z)+1$,
we may assume $d_h(z,x)\leq d_h(y,x)+1$. Now $(x|z)\ge (t-1.5+t-1.5-(d_h(y,x)+1))/2=
((t+1.5)+(t+1.5)-d_h(y,x))/2-3.5\ge (x|y)-4$.
\par Arrange points $x$, $y$, and $z$ as $u$, $v$, and $w$
so that $s=d_h(u,x_0)\leq m=d_h(v,x_0)\leq l=d_h(w,x_0)$.
We may assume $s\leq l-4$ (otherwise $s$, $m$, and $l$ lie within $1.5$ from
$(s+l)/2$). 
 
 Case 1: $m\leq l-2$. Now $d_h(w,v)$ must be within $1$ from $l$, so
 $2(w|v)= m+l-d_h(w,v)$ is contained between $m-1$ and $m+1$.
 Similarly, $s-1\leq 2(u|w)\leq s+1$.
 Since $d_h(u,v)\leq m+1$, $2(u|v)\ge m+s-(m+1)=s-1$.
 The only possibility for the smallest of Gromov products for the triple $uvw$
 to be less than the medium one minus $4$ is if $s\leq m-2$.
 In that case $d_h(u,v)\ge m-1$, so  $s+1=m+s-(m-1)\ge 2(u|v)\ge m+s-(m+1)=s-1$
 and the smallest of Gromov products for the triple $uvw$ is
larger than the medium one minus $4$.
 
 Case 2: $m\ge l-1$. Now $d_h(u,v)$ must be within $1$ from $m$, so
 $2(u|v)= m+s-d_h(u,v)$ is contained between $s-1$ and $s+1$.
 Similarly, $s-1\leq 2(u|w)\leq s+1$.
 Since $d_h(w,v)\leq l+1$, $2(w|v)\ge m+l-(l+1)=m-1\ge s+2$, and the smallest of Gromov products for the triple $uvw$ is
larger than the medium one minus $4$.
\par
To prove $\dim_{NA}(X,d_h)\leq n$ we plan to define covers $\VV_r$ of $(X,d_h)$ such that
$\mesh(\VV_r)\leq r$, $L(\VV_r)\ge r/4$, and $m(\VV_r)\leq n+1$.
If $r\leq 4$, we define $\VV_r$ as all singletons of $X$,
otherwise we put $\VV_r=\UU_i$ with $i$ being the integral part of $r$.
Indeed, for $r > 4$, $\mesh(\UU_i)\leq i\leq r$ and $L(\UU_i)\ge i-1\ge r-2\ge r/4$.

\par
Let us show the Gromov boundary $\partial_\infty X$
of $X$ consists of one point. Given two sequences of points $\{x_i\}$ and $\{y_i\}$ such that
$(x_m|x_k)\to \infty$ and $(y_m|y_k)\to\infty$ we need to prove $(x_m,y_m)\to\infty$. 
However, $d_h(x_m,y_m)$ is smaller than $\max(d_h(x_m,x_0),d_h(y_m,x_0))+2$,
so $2(x_m|y_m)\ge \min(d_h(x_m,x_0),d_h(y_m,x_0))-2\to\infty$.
\par
Since b)$\implies$c) and c)$\implies$a) are obvious, we are done.
\hfill $\blacksquare$

\end{document}